\documentclass[review,nonatbib]{elsarticle}
\usepackage{amsmath}
\usepackage{amssymb}

\makeatletter
\let\c@author\undefined
\makeatother
\usepackage[backend=biber]{biblatex}
\usepackage{tikz}
\usetikzlibrary{positioning}
\usepackage{subfiles}

\newcommand{\R}{\mathbb{R}}

\renewcommand{\det}{\operatorname{det}}
\newcommand{\vol}{\operatorname{vol}}

\renewcommand{\top}{\mathsf{T}}

\newcommand{\RR}{\mathbb{R}}

\let\oldbigwedge\bigwedge
\renewcommand{\bigwedge}{\mathop{\oldbigwedge}\nolimits}

\newtheorem{lemma}{Lemma}
\newtheorem{theorem}{Theorem}
\newtheorem{corollary}{Corollary}
\newtheorem{proposition}{Proposition}
\newtheorem{definition}{Definition}

\newproof{proof}{Proof}

\begin{document}

\begin{frontmatter}
\title{Determinant-Based Error Bounds for CUR Matrix Approximation: Oversampling and Volume Sampling}

\author[3]{Frank de Hoog\corref{cor1}}
\ead{Frank.deHoog@csiro.au}

\author[4]{Markus Hegland}
\ead{markus.hegland@anu.edu.au}

\cortext[cor1]{Corresponding author}

\affiliation[3]{organization={CSIRO Data61},
                city={Canberra},
                postcode={2601}, 
                country={Australia}}

\affiliation[4]{organization={The Australian National University},
                city={Canberra},
                postcode={0200}, 
                country={Australia}}

\begin{abstract}
We derive error bounds for CUR matrix approximation using determinant-based methods that relate local projection errors to global approximation quality. For general matrices, we establish determinant identities for bordered Gramian matrices that decompose CUR approximation errors into interpretable local components. These identities connect projection errors onto submatrix column spaces directly to determinants, providing geometric insight into the degradation of the approximation. We develop a probabilistic framework based on volume sampling that yields interpolation-type error bounds quantifying the benefits of oversampling: when $r > k$ rows are sampled to construct a rank-$k$ approximation, the expected error factor decreases linearly from $(k+1)^2$ (no oversampling) to $(k+1)$ (full oversampling). Our analysis establishes that the expected squared error is bounded by this interpolation factor times the squared error of the best rank-$k$ approximation, directly connecting the quality of the CUR approximation to that of the optimal low-rank approximation. The framework applies to both CUR decomposition for general matrices and the Nyström method for symmetric positive semi-definite matrices, providing a unified theoretical foundation for the analysis of determinant-based low-rank approximation.
\end{abstract}

\begin{keyword}
Matrix approximation \sep CUR decomposition \sep determinants \sep volume sampling \sep oversampling
\end{keyword}

\end{frontmatter}

\maketitle

\tableofcontents

\section{Introduction}\label{sec:1}

\subsection{Low-Rank Matrix Approximation and the CUR Framework}

Low-rank matrix approximation is a fundamental computational primitive in modern data science, underpinning diverse applications ranging from recommendation systems and image compression \cite{allen2024} to kernel methods via the Nyström approach \cite{WilS00}, hierarchical solvers for PDEs \cite{bebendorf2008}, reduced-order modelling in computational engineering~\cite{benner2015survey}, data compression in network monitoring~\cite{wang2024efficient}, and large-scale machine learning \cite{mahoney2011randomized}.

Given a large matrix $M \in \mathbb{R}^{m \times n}$, the central goal is to find a \emph{rank-$k$ approximation} with $k \ll \min(m,n)$ that captures the essential structure of the data while dramatically reducing storage and computational requirements.

The classical solution is the \emph{truncated singular value decomposition (SVD)}, which provides the optimal rank-$k$ approximation in all Schatten $p$-norms---a fundamental result established by Eckart and Young \cite{eckart1936}. This optimality has made the truncated SVD the gold standard for low-rank approximation \cite{golub2013matrix}. However, computing truncated SVDs requires access to the entire matrix and incurs computational costs that become prohibitive for massive datasets. Moreover, the resulting singular vectors are abstract linear combinations of the original data elements, making them difficult to interpret in many practical applications \cite{mahoney2009cur}.

An alternative paradigm seeks to approximate matrices directly from subsets of their rows and columns. In its simplest form, one selects $k$ rows and $k$ columns to form an intersection submatrix $A = M_{I,J}$. If $A$ has full rank, the optimal reconstruction of this type takes the form
\[M \approx C U R, \quad C = M_{:,J}, \; R = M_{I,:}, \; U = A^{-1},\]

where $A^{-1}$ denotes the inverse of $A$. This construction, known as a \emph{skeleton} or \emph{pseudo-skeleton approximation} \cite{goreinov1997pseudoskeleton}, has the crucial property that it interpolates the exact entries of $M$ in the chosen rows and columns.

In the more general \emph{CUR decomposition} with \emph{oversampling}, more than $k$ rows or columns may be selected (i.e., $r > k$) to improve robustness \cite{park2024accuracy}. The theoretically optimal choice for the middle factor is
\begin{equation} U^\star = C^{+} M R^{+},\end{equation}

as demonstrated by Stewart \cite{stewart1999} and employed in the general CUR framework \cite{mahoney2009cur}. This choice minimizes the Frobenius norm error over all possible $U$. However, computing $U^\star$ requires access to the entire matrix $M$, which defeats the purpose of the submatrix approach. A widely used and computationally efficient alternative retains the simple choice $U = A^{+}$ (the pseudoinverse of the rectangular submatrix), for which sharp deterministic error bounds can still be established. This fundamental trade-off between the optimal but expensive choice $U^\star$ and the computationally tractable $A^+$ is central to the theoretical development and practical algorithms studied in this paper, particularly in the oversampling regime where $r > k$.

\subsection{Determinant Identities and Volume Sampling}

The principal contribution of this paper is a comprehensive analysis of CUR and Nyström approximation errors in the Frobenius norm, achieved by relating global error bounds to local quantities expressed through determinants of Gramian matrices. Our approach builds on the observation that determinants naturally encode geometric information about matrix approximations. For instance, when augmenting a submatrix $A$ with an additional column $b$, we have the fundamental identity~\cite{horn2012matrix}

\[\det([A,b]^\top [A,b]) = \det(A^\top A)\,\|(I-AA^+)b\|^2.\]

This identity reveals that local determinants directly encode projection errors---the right-hand side explicitly separates the volume contribution of the existing subspace from the residual component orthogonal to it. Similar formulas hold when adding rows, or simultaneously adding both a row and a column, and these can be unified using compound matrices and the Cauchy--Binet theorem. These local determinant identities form the fundamental building blocks for the global error analysis developed in subsequent sections. The analysis is termed 'local' because these determinant identities capture the approximation errors introduced by individual rows and columns from the original matrix, allowing us to understand how each additional data element affects the overall CUR decomposition quality.

A key insight emerging from this determinant-based analysis is a family of \emph{interpolation-type error bounds} that precisely quantify how oversampling improves approximation quality. When $r = k$ (no oversampling), the expected Frobenius error satisfies a $(k+1)^2$-factor bound relative to the best rank-$k$ approximation, recovering the classical result of \cite{goreinov1997pseudo}. As the oversampling parameter $r$ increases, this factor decreases linearly, reaching a $(k+1)$-factor when $r = m$. This demonstrates that oversampling rows and columns not only enhances numerical stability but also yields provably sharper error bounds, consistent with recent deterministic analyses \cite{cortinovis2020} and randomized methods \cite{halko2011}. While the bounds on the Frobenius norm of the CUR error also apply to the Nyström method, they are consistent with bounds previously obtained from the eigenvalue decomposition \cite{GittensMahoney2016}.

The analysis relies on connections between matrix approximation errors and elementary symmetric polynomials in the singular values, established through compound matrix theory (detailed in Section 3).
\subsection{Historical Development and Recent Advances}

An extensive review of low-rank matrix approximations including applications and algorithms can be found in \cite{KisS17}.

The theoretical foundation for CUR and related decompositions has evolved through several distinct phases. The earliest roots lie in classical projection theory, where orthogonal projections $AA^+$ are used to approximate matrices of the form $[A,B]$, with error bounds expressed through principal angles \cite{davis1970rotation}. These results were subsequently popularized in the numerical linear algebra community through the influential textbook of Golub and Van Loan \cite{golub2013matrix}.

A significant advance came with the introduction of rank-revealing factorizations, particularly RRQR decompositions, which provide deterministic methods for selecting columns that preserve numerical rank and approximate singular values. Key contributions in this area include the work of Gu and Eisenstat~\cite{gu1996efficient} and other foundational developments in numerical linear algebra. This development strongly influenced later column subset selection strategies and established connections to condition number bounds.

Concurrently, the emergence of randomized numerical linear algebra brought new perspectives to matrix approximation. Notable contributions include the CUR algorithms of Drineas, Mahoney, and collaborators \cite{drineas2008relative, mahoney2009cur} and the comprehensive survey by Halko, Martinsson, and Tropp \cite{halko2011}. These randomized approaches achieve relative-error guarantees in expectation by sampling according to leverage scores or related probability distributions, offering computational advantages for very large matrices. An error bound for the projection method ($A=C$) is established in~\cite{deshpande2006matrix}
where it is shown that there exists a submatrix $A\in\R^{m,k}$ of $M$ satisfying 
\[\|M - AA^+M\|_F^2 \le (k+1)\,\|M - M_k\|_F^2,\]
where $M_k$ denotes the best rank-$k$ approximation of $M$. This established the $(k+1)$-factor that continues to guide theoretical developments in the field.

Recent developments have seen the rise of deterministic algorithms that achieve comparable error bounds without explicit volume computations, most prominently in the work of Cortinovis and Kressner \cite{cortinovis2020}. Simultaneously, there has been renewed interest in determinantal point processes and volume sampling as principled approaches to row and column selection \cite{derezinski2021determinantal}, with connections to fast sampling algorithms and improved approximation guarantees.

Modern algorithmic developments include efficient approximate sampling methods for determinantal point processes \cite{grosse2024greedy} and scalable implementations for nonsymmetric cases \cite{han2022scalable}. Additionally, hybrid approaches that combine deterministic and randomized strategies have emerged to balance computational efficiency with approximation accuracy \cite{abdelgawad2025hycur}.

\subsection{Contributions and Paper Organization}

This paper provides a unified theoretical framework that bridges local determinant identities with global approximation bounds for both CUR and Nyström methods. Our main contributions are as follows:

\textbf{Local Error Analysis:} In Section~2, we develop a comprehensive local analysis based on determinant identities for bordered Gramian matrices. By augmenting a chosen submatrix with additional rows and/or columns, we derive explicit formulas that decompose global approximation errors into interpretable local components. Our central result (Proposition~3) establishes that for a bordered matrix

\[X = \begin{bmatrix} A & b \\ c^\top & d \end{bmatrix},\]

one has the explicit identity

\[\det(X^\top X) = \det(A^\top A + cc^\top)\,\|(I - AA^+)b\|^2 + \det(A^\top A)\,(d - c^\top A^+ b)^2.\]

This formula directly relates $\det(X^\top X)$ to squared local error terms, providing geometric insight into how matrix approximations degrade when additional data is incorporated.

\textbf{Local Deterministic Errors:} We establish sharp deterministic bounds that require only average-volume rather than maximal-volume assumptions. Specifically, Proposition~7 shows that

\[\|X_{\mathrm{CUR}} - X\|^2 \le \frac{(r+1)(k+1)}{r+1-k}\,\lambda_{\min}(X^\top X),\]

provided the selected submatrix has squared volume at least as large as the average over all possible submatrices of the same size. This represents a significant relaxation of traditional maximal-volume requirements.

\textbf{Global Probabilistic Analysis:} Section~3 develops a probabilistic framework based on volume sampling that connects local determinant identities to global approximation guarantees. We derive interpolation-type error bounds that precisely quantify the benefits of oversampling: the error factor transitions smoothly from $(k+1)^2$ when $r = k$ to $(k+1)$ when $r = m$, providing both theoretical insight and practical guidance for algorithm design.

In particular, we establish bounds on the expected squared Frobenius error directly in terms of the singular values of $M$, providing a direct connection between the CUR approximation error and the singular value tail:
\[\mathbb{E}\left(\|M-M_{\text{CUR}}\|_F^2\right) \leq \left(\frac{m-r}{m-k}(k+1)^2 + \frac{r-k}{m-k}(k+1) \right) \sum_{i=k+1}^n \sigma_i^2.\]

Together, these contributions provide a comprehensive understanding of the mechanisms governing CUR and Nyström approximation quality, clarifying the role of determinants, volume sampling, and oversampling in achieving accurate low-rank approximations.

The remainder of the paper is organised as follows. Section~2 presents the local determinant analysis and establishes connections to compound matrices via the Cauchy--Binet theorem. Section~3 develops the global volume sampling framework and derives the main approximation bounds. We conclude with a discussion of the implications for algorithm design and directions for future research.

\ifdefined\standalonecompile
  \printbibliography
\fi

\setcounter{section}{1}
\section{Local error analysis for CUR}
\label{sec:local}
This section develops a local (bordered) determinant-based analysis of the CUR approximation. 
 The core idea is simple and geometric: starting from a base submatrix
\begin{equation}\label{eq:A}
A = M_{I,J},\qquad A\in\RR^{r\times k},\; r\ge k,
\end{equation}
we study the effect of adding a single column, a single row, or both. The corresponding extended matrices
\begin{equation}\label{eq:XY}
Y = [A\; b],\qquad Z = \begin{bmatrix}A \\ c^\top\end{bmatrix},\qquad X=\begin{bmatrix}A & b\\ c^\top & d\end{bmatrix}
\end{equation}
encode the local geometry: their Gram determinants measure squared volumes (compound norms) and, as we show below, factor precisely into a product of an "ambient" determinant and the squared norm of a residual term. These identities are the algebraic heart of our local error estimates and feed directly into the bounds in Section~3.

Throughout this section we use \(A,b,c,d\) to denote local bordered blocks of the
matrix \(M\), while \(I,J\) index the corresponding selected rows and columns.
The bordered notation keeps the local algebra transparent, the indexed notation reappears in Section~3 when these local results are lifted to the global setting.

We organise the section as follows. Section~2.1 establishes determinant identities for bordered Gramians. Section~2.2 recasts these identities using compound matrices and the Cauchy--Binet theorem to expose the geometric interpretation. Section~2.3 converts the identities into deterministic CUR error bounds.

\subsection{Determinant identities for bordered Gramians}

We begin with two elementary but fundamental identities (column and row additions). Throughout we assume the usual full-rank conditions so that inverses or pseudoinverses are well defined; the formulae extend by continuity to rank-deficient cases and we mention the appropriate adjugate corrections where useful.

\begin{proposition}[Adding a column]\label{prop:add-col}

Let \(A \in \mathbb{R}^{r\times k}\) with \(r \ge k\) and full column rank \(k\), and let \(b \in \mathbb{R}^r\).
 Set $Y=[A\; b] \in\RR^{r\times(k+1)}$. Then
\begin{equation}\label{eq:add-col}
\det\bigl(Y^\top Y\bigr) 
= \det\bigl(A^\top A\bigr)\,\lVert (I-AA^+)b\rVert^2.
\end{equation}
\end{proposition}

The determinant of the Gramian of $Y$ thus factors into the determinant of the Gramian of $A$ times the squared norm of the projection of $b$ onto the orthogonal complement of the range of $A$.

\begin{proof}
Write the bordered Gramian in block form:
$$Y^\top Y=\begin{bmatrix}A^\top A & A^\top b\\ b^\top A & b^\top b\end{bmatrix}.$$
Applying the Schur complement (Cauchy expansion) identity gives
$$\det(Y^\top Y)=\det(A^\top A)\bigl(b^\top b - b^\top A(A^\top A)^{-1}A^\top b\bigr).$$
Observe that $A(A^\top A)^{-1}A^\top=AA^+$, and rearranging yields \eqref{eq:add-col}.
\end{proof}

\begin{proposition}[Adding a row]\label{prop:add-row}
Let $A\in\RR^{r\times k}$ have full column rank $k$ and $c\in\RR^k$. Define $Z=\begin{bmatrix}A\\ c^\top\end{bmatrix}\in\RR^{(r+1)\times k}$. Then
\begin{equation}\label{eq:add-row}
\det\bigl(Z^\top Z\bigr)=\det\bigl(A^\top A + c c^\top\bigr)=\det(A^\top A)\bigl(1 + c^\top(A^\top A)^{-1}c\bigr).
\end{equation}
\end{proposition}

\begin{proof}
Direct computation gives $Z^\top Z=A^\top A + c c^\top$; the displayed identity follows from the matrix determinant lemma (rank-one update formula)~\cite[Thm.~0.7.2]{horn2012matrix}.
\end{proof}

\noindent\textbf{Remark.} If $A$ is rank-deficient the right-hand side must be interpreted using the adjugate; the rank-one update identity continues to hold in an appropriately generalised form..

When both a row and a column are added, one obtains a two-term decomposition: one term from the column residual and another from a scalar-valued Schur complement.

\begin{proposition}[Adding a row and a column]\label{prop:add-both}
Let $A\in\RR^{r\times k}$ have full column rank $k$, and let $b\in\RR^r$, $c\in\RR^k$, $d\in\RR$. Set
$$X=\begin{bmatrix}A & b\\ c^\top & d\end{bmatrix}\in\RR^{(r+1)\times(k+1)}.$$ 
Denote $u=(I-AA^+)b$ and $\gamma=d-c^\top A^+ b$. Then
\begin{equation}\label{eq:add-both}
\det\bigl(X^\top X\bigr) = \det\bigl(A^\top A + cc^\top\bigr)\,\lVert u\rVert^2 + \det\bigl(A^\top A\bigr)\,\gamma^2.
\end{equation}
\end{proposition}

The squared volume of $X$ splits into contributions from the column residual $u$ and the scalar Schur complement $\gamma$.

\begin{proof}
Write $b=AA^+b+u$ and $d=c^\top A^+ b + \gamma$, and factor
$$X=\begin{bmatrix}A & u\\ c^\top & \gamma\end{bmatrix}\begin{bmatrix}I & A^+b\\ 0 & 1\end{bmatrix}=:WR.$$ 

Since \(\det(R) = 1\), we have \(\det(X^{\top}X) = \det(W^{\top}W)\).
Writing $B=A^\top A + cc^\top$ one checks
$$W^\top W=\begin{bmatrix}B & \gamma c\\ \gamma c^\top & \lVert u\rVert^2+\gamma^2\end{bmatrix},$$
and applying the bordered determinant formula yields the stated decomposition after a short algebraic simplification using $B-cc^\top=A^\top A$.
\end{proof}

\subsection{Compound matrices and a geometric interpretation}

For a matrix $A \in \RR^{r \times k}$ with $r \geq k$, the \emph{$k$-th order compound matrix} 
$C_k(A)$ is the $\binom{r}{k} \times 1$ column vector whose entries are all $k$-th order minors of $A$. 
Equivalently, $C_k(A)$ is the exterior product (wedge product) $A_{:,1} \wedge \cdots \wedge A_{:,k}$ of the columns of $A$. Similarly, for the extended matrices $Y \in \RR^{r \times (k+1)}$ and $X \in \RR^{(r+1) \times (k+1)}$, we have $C_{k+1}(Y) \in \RR^{\binom{r}{k+1} \times 1}$ and $C_{k+1}(X) \in \RR^{\binom{r+1}{k+1} \times 1}$, again column vectors of the corresponding minors.

The Cauchy--Binet theorem shows that Gram determinants equal squared norms of these compound matrices.

\begin{proposition}[Cauchy--Binet identities]\label{prop:cauchy-binet}
Let $A\in\RR^{r\times k}$, $Y\in\RR^{r\times(k+1)}$, $X\in\RR^{(r+1)\times(k+1)}$ with $r\ge k$. Then
\begin{align*}
\det(A^\top A)&=\lVert C_k(A)\rVert^2,\\
\det(Y^\top Y)&=\lVert C_{k+1}(Y)\rVert^2,\\
\det(X^\top X)&=\lVert C_{k+1}(X)\rVert^2.
\end{align*}
\end{proposition}

The identities in Propositions~\ref{prop:add-col}--\ref{prop:add-both} now admit an immediate interpretation: the compound matrix of the extended configuration equals the wedge product of the base compound matrix with a compound matrix encoding the local residual, and the determinant factorisation is the squared-norm identity in the exterior algebra.

As a corollary we isolate the increment in Gram determinant when adding a row:
\begin{lemma}\label{lem:gram-increment}
Let $X\in\RR^{(r+1)\times(k+1)}$ and $Y=X_{1:r,:}$. If $r\ge k+1$ then
\begin{equation}\label{eq:gram-inc}
\det(X^\top X)-\det(Y^\top Y)=\sum_{H\in\mathcal{S}_k(r+1)}\mathbf{1}_{H\not\subset[r]}\det(X_{H,:})^2.
\end{equation}
\end{lemma}

\begin{proof}
Immediate from Cauchy--Binet: separate minors involving the new row from those that do not.
\end{proof}

\subsection{Local CUR error bounds}

We now connect the determinant identities to local CUR errors. Consider a matrix $X\in\RR^{(r+1)\times(k+1)}$ and suppose we choose an $r\times k$ submatrix $A_{i^*,j^*}$ (obtained by removing one row and one column) satisfying a natural averaging criterion. The following propositions show how this choice controls the local error in terms of the smallest eigenvalue of $X^\top X$. These bounds are not used in the derivation of the expected errors under volume sampling in Section~3, and the subsection may therefore be skipped by readers interested only in the global error bounds. Nonetheless, these bounds demonstrate that strong deterministic error bounds for local errors can be obtained.

We begin with a foundational eigenvalue bound for positive definite matrices.

\begin{proposition}\label{prop:eig-bound-square}
Let $G \in \RR^{(k+1) \times (k+1)}$ be positive definite. For each $i \in \{1,\ldots,k+1\}$, let $F_i$ denote the principal $k \times k$ submatrix formed by deleting row $i$ and column $i$ of $G$. If an index $i^*$ satisfies
$$\det(F_{i^*}) \geq \frac{1}{k+1}\sum_{i=1}^{k+1} \det(F_i),$$
then
\begin{equation}\label{eq:eig-bound-square}
\frac{\det(G)}{\det(F_{i^*})} \leq (k+1) \cdot \lambda_{\min}(G).
\end{equation}
\end{proposition}

\begin{proof}
By hypothesis, 
$$\det(F_{i^*}) \geq \frac{1}{k+1}\sum_{i=1}^{k+1} \det(F_i),$$
which implies
$$\frac{\det(G)}{\det(F_{i^*})} \leq \frac{(k+1)\det(G)}{\sum_{i=1}^{k+1} \det(F_i)}.$$
The theory of elementary symmetric polynomials gives
$$\sum_{i=1}^{k+1} \det(F_i) = s_k(\lambda_1,\ldots,\lambda_{k+1}) = \det(G)\sum_{j=1}^{k+1}\frac{1}{\lambda_j},$$
where $s_k$ is the $k$-th elementary symmetric polynomial in the eigenvalues $\lambda_j$ of $G$. Therefore
$$\frac{(k+1)\det(G)}{\sum_{i=1}^{k+1} \det(F_i)} = \frac{k+1}{\sum_{j=1}^{k+1}\frac{1}{\lambda_j}}.$$
Since $\sum_{j=1}^{k+1}\frac{1}{\lambda_j} \geq \frac{1}{\lambda_{\min}(G)}$ (the sum of positive reciprocals is at least the reciprocal of the minimum), we obtain the claimed bound.
\end{proof}

The next proposition extends this to rectangular matrices via a multiplicity counting argument.

\begin{proposition}\label{prop:rect-bound}
Let $X \in \RR^{(r+1) \times (k+1)}$ with $r \geq k$. For $j \in \{1,\ldots,k+1\}$, let $Y_j \in \RR^{(r+1) \times k}$ be $X$ with column $j$ removed. For $i \in \{1,\ldots,r+1\}$ and $j \in \{1,\ldots,k+1\}$, let $A_{i,j} \in \RR^{r \times k}$ be $Y_j$ with row $i$ removed. If indices $(i^*,j^*)$ satisfy
$$\det(A_{i^*,j^*}^\top A_{i^*,j^*}) \geq \frac{1}{(r+1)(k+1)}\sum_{i,j} \det(A_{i,j}^\top A_{i,j}),$$
then
\begin{equation}\label{eq:rect-bound}
\frac{\det(X^\top X)}{\det(A_{i^*,j^*}^\top A_{i^*,j^*})} \leq \frac{(r+1)(k+1)}{r+1-k} \cdot \lambda_{\min}(X^\top X).
\end{equation}
\end{proposition}

\begin{proof}
By hypothesis,
$$\det(A_{i^*,j^*}^\top A_{i^*,j^*}) \geq \frac{1}{(r+1)(k+1)}\sum_{i,j}\det(A_{i,j}^\top A_{i,j}),$$
which yields
$$\frac{\det(X^\top X)}{\det(A_{i^*,j^*}^\top A_{i^*,j^*})} \leq \frac{(r+1)(k+1)\det(X^\top X)}{\sum_{i,j}\det(A_{i,j}^\top A_{i,j})}.$$
We now apply Cauchy--Binet to relate the sums of determinants. Write 
$$\det(Y_j^\top Y_j) = \sum_{I \in \mathcal{S}_k(r+1)} \det(W_{I,j})^2,$$ 
where $W_{I,j}$ are $k \times k$ submatrices of $Y_j$ selecting $k$ rows from the $(r+1)$ available. 
Here \(S_k(r+1)\) denotes the set of all \(k\)-element subsets of \(\{1,\dots,r+1\}\),
corresponding to row selections in the Cauchy–Binet expansion.
Thus the above sum has $\binom{r+1}{k}$ terms. Similarly,
$$\det(A_{i,j}^\top A_{i,j}) = \sum_{I \in \mathcal{S}_k([r+1]\setminus\{i\})} \det(W_{I,j})^2$$
sums over $\binom{r}{k}$ terms. The key observation is that each $k \times k$ submatrix $W_{I,j}$ appears \emph{once} in the $Y_j$ sum and \emph{$(r+1-k)$ times} in the $A_{i,j}$ sums (once for each row $i \notin I$). Summing over all $(i,j)$ pairs gives
$$\sum_{i,j}\det(A_{i,j}^\top A_{i,j}) = (r+1-k)\sum_j\det(Y_j^\top Y_j).$$
Therefore
$$\frac{(r+1)(k+1)\det(X^\top X)}{\sum_{i,j}\det(A_{i,j}^\top A_{i,j})} = \frac{(r+1)(k+1)}{r+1-k} \cdot \frac{\det(X^\top X)}{\sum_j\det(Y_j^\top Y_j)}.$$
Applying Proposition~\ref{prop:eig-bound-square} to the Gramian $X^\top X$ and its principal submatrices $Y_j^\top Y_j$ completes the proof.
\end{proof}

The preceding bounds translate directly into a local CUR error estimate.

\begin{proposition}\label{prop:local-cur-final}
Let $X \in \RR^{(r+1) \times (k+1)}$ with $r \geq k$, and let $X_{\mathrm{CUR}}$ be the CUR approximation of $X$ formed using a submatrix $A_{i^*,j^*}$ satisfying
$$\det(A_{i^*,j^*}^\top A_{i^*,j^*}) \geq \frac{1}{(r+1)(k+1)}\sum_{i,j}\det(A_{i,j}^\top A_{i,j}).$$
Then
\begin{equation}\label{eq:local-cur-final}
\|X_{\mathrm{CUR}} - X\|_F^2 \leq \frac{(r+1)(k+1)}{r+1-k} \cdot \lambda_{\min}(X^\top X).
\end{equation}
\end{proposition}

\begin{proof}
From Proposition~\ref{prop:add-both}, the determinant relation for the bordered Gramian gives
$$\det(X^\top X) = \det(A_{i^*,j^*}^\top A_{i^*,j^*} + cc^\top)\lVert u\rVert^2 + \det(A_{i^*,j^*}^\top A_{i^*,j^*})\gamma^2,$$
where $u = (I - A_{i^*,j^*}A_{i^*,j^*}^+)b$ and $\gamma = d - c^\top A_{i^*,j^*}^+ b$ encode the local CUR residuals. This implies
$$\|X_{\mathrm{CUR}} - X\|_F^2 \leq \frac{\det(X^\top X)}{\det(A_{i^*,j^*}^\top A_{i^*,j^*})}.$$
Applying Proposition~\ref{prop:rect-bound} yields the claimed bound.
\end{proof}

\noindent\textbf{Remarks.}
\begin{itemize}
\item The factor $(r+1)(k+1)/(r+1-k)$ quantifies the combinatorial inflation from removing rows and columns and is mild when $k \ll r$.
\item The bound is constructive: selecting a submatrix whose Gram determinant is at least the average of all such submatrices guarantees a controlled error tied to $\lambda_{\min}(X^\top X)$.
\item In the square invertible case ($r=k$, $A$ invertible), the $u$ term vanishes and the bound reduces to a classical Schur complement estimate.
\end{itemize}

\paragraph{Additional remark.}
The bounds established in this subsection apply equally to the symmetric positive definite case.  When \(X\) is symmetric and positive definite, the Gramian \(X^{\top}X\) coincides with \(X^2\), and the intermediate lemma used in Proposition~\ref{prop:rect-bound} already covers this situation.  It follows that the resulting determinant ratio and Frobenius-norm error bounds coincide with those obtained by analysing the symmetric CUR (Nyström) approximation directly.  Hence the present framework unifies both the general rectangular and the symmetric positive definite settings under the same local determinant identities.

\medskip
\noindent\textbf{Connection to the literature.} The use of orthogonal projections $AA^+$ to approximate matrices of the form $\begin{bmatrix}A & B\end{bmatrix}$ has deep roots in classical least squares theory and was popularised in numerical linear algebra through the textbook by Golub and Van Loan~\cite{golub2013matrix}. The projection error $\lVert b - AA^+b\rVert$ is geometrically related to the angle $\theta$ between $b$ and the column space of $A$ via $\lVert b - AA^+b\rVert \leq \sin(\theta)\lVert b\rVert$, cf.~Davis and Kahan~\cite{davis1970rotation}.

Bounds on projection error in terms of leverage scores emerged in the randomised numerical linear algebra literature, particularly through work by Drineas, Mahoney, and collaborators~\cite{drineas2008relative}. The first rigorous proof that a $k$-column submatrix of maximal volume provides the deterministic bound $\lVert M - AA^+M\rVert_F^2 \leq (k+1)\lVert M - M_k\rVert_F^2$ was given by Deshpande et al.~\cite{deshpande2006matrix}.

\medskip
\noindent\textbf{Summary.} This section established exact algebraic identities linking determinants of bordered Gramians to projection residuals and scalar Schur complements, and translated those identities into interpretable local CUR error bounds. The determinant identities (Propositions~\ref{prop:add-col}--\ref{prop:add-both}) provide the algebraic foundations; the Cauchy--Binet reformulation clarifies the geometric meaning in terms of compound matrices; and the combinatorial counting arguments convert local identities into deterministic inequalities. These local building blocks are used in Section~3 to derive global expected-error bounds under volume sampling and related selections.

\ifdefined\standalonecompile
  \printbibliography
\fi

\setcounter{section}{2}

\section{Main results for general matrices}
\label{sec:results}
\sloppy
This section derives explicit expressions and upper bounds for the error in the CUR approximation of a matrix $M$.
We introduce the probability space and random variables underlying the randomised CUR model, and derive the expected squared Frobenius norm of the error of the CUR approximation. Throughout this section we assume $M \in \mathbb{R}^{m \times n}$ with $m \ge n$.

We now extend the bordered determinant identities of Section~2 to the global matrix setting.
The same determinant and compound-matrix framework applies, but here the index sets are random and the results are expressed as expectations over the sampling distribution.

\subsection{Volume sampling and CUR approximation}

Let $M \in \mathbb{R}^{m \times n}$ be a fixed matrix, $1 \leq k \leq r \leq m$ and $k \leq n$. We define a sampling distribution over index sets and use it to construct a randomised CUR approximation of $M$. In addition, let $\mathcal{S}_r(m)$ denote the set of all subsets of $[m]$ of size $r$, and $\mathcal{S}_k(n)$ the set of all subsets of $[n]$ of size $k$.

Throughout this section we use the same notation as in Section~2; the $k$-th compound matrix of $M$ is denoted $C_k(M)$, and the index sets are written $\mathcal{S}_r(m)$ and $\mathcal{S}_k(n)$ for subsets of $\{1,\dots,m\}$ and $\{1,\dots,n\}$, respectively.

\begin{definition}[Volume sampling probability]
    \label{def:volume-sampling}
  For each pair $(I, J) \in \mathcal{S}_r(m) \times \mathcal{S}_k(n)$, let the \emph{volume sampling probability} be
$$
    p(I, J) = \zeta^{-1} \det(M_{I,J}^\top M_{I,J}), \quad (I, J) \in \mathcal{S}_r(m) \times \mathcal{S}_k(n),
  $$
  where the \emph{volume sampling normalisation factor} is
  $$
    \zeta = \sum_{(I,J) \in \mathcal{S}_r(m) \times \mathcal{S}_k(n)} \det(M_{I,J}^\top M_{I,J}).
  $$
\end{definition}

\medskip
The next result provides a closed-form expression for the normalisation factor, analogous to the local determinant identities of Section~2.2.
\begin{theorem}[Volume sampling normalisation factor]\label{thm:volume-sampling-normalisation}
  Let $M, m, r, k$ be as in Definition~\ref{def:volume-sampling}. 
  Then the volume-sampling normalisation factor $\zeta$ satisfies
  $$\zeta = \binom{m-k}{r-k}\,\|C_k(M)\|_F^2.$$
\end{theorem}

\begin{proof}
From the definition,
$$
\zeta = \sum_{(I,J) \in \mathcal{S}_r(m) \times \mathcal{S}_k(n)} \vol(M_{I,J})^2.
$$
Applying the Cauchy--Binet formula gives
$$
\vol(M_{I,J})^2 = \sum_{K \in \mathcal{S}_k(m)} \mathbf{1}_{K \subseteq I} \vol(M_{K,J})^2.
$$
Hence 
$$
\zeta =  \sum_{I \in \mathcal{S}_r(m)} \sum_{J \in \mathcal{S}_k(n)}
  \sum_{K \in \mathcal{S}_k(m)} \mathbf{1}_{K \subseteq I} \vol(M_{K,J})^2.
$$

Reordering the summations yields
$$
\zeta = \sum_{J \in \mathcal{S}_k(n)} \sum_{K \in \mathcal{S}_k(m)} \vol(M_{K,J})^2 \left( \sum_{I \in \mathcal{S}_r(m)} \mathbf{1}_{K \subseteq I} \right).
$$

The inner sum 
$$
\sum_{I \in \mathcal{S}_r(m)} \mathbf{1}_{K \subseteq I}
$$
counts the number of sets $I$ of size $r$ containing a fixed set $K$ of size $k$.
Any set $I$ of size $r\geq k$ containing $K$ is formed by including $r-k$ elements of the set $[m]\setminus K$ to $K$. As the cardinality of the $[m]\setminus K$ is $m-k$ the number of resulting sets $I$ is
$$
\binom{m-k}{r-k}.
$$

Thus,
$$
\zeta = \binom{m-k}{r-k} \sum_{J \in \mathcal{S}_k(n)} \sum_{K \in \mathcal{S}_k(m)} \vol(M_{K,J})^2.
$$

Finally, by the definition of the Frobenius norm of the $k$-th compound matrix,
$$
\| C_k(M) \|_F^2 = \sum_{K \in \mathcal{S}_k(m)} \sum_{J \in \mathcal{S}_k(n)} \vol(M_{K,J})^2,
$$
and therefore
$$
\zeta = \binom{m-k}{r-k} \| C_k(M) \|_F^2.
$$
\end{proof}

This shows that the volume-sampling weight is proportional to the squared Frobenius norm of the $k$-th compound matrix.

\begin{definition}[Random CUR approximation]
Let $I \in \mathcal{S}_r(m)$, $J \in \mathcal{S}_k(n)$ be selected using a sampling procedure (e.g., volume sampling).
    
Once $I$ and $J$ are chosen, the corresponding CUR approximation of $M$ used here is defined as
$$\mathsf{CUR} = M_{:,J} M_{I,J}^+ M_{I,:}.$$
This approximation depends on $I$ and $J$ and therefore varies from sample to sample.
\end{definition}

Expectations are taken over this sampling procedure. The expected squared Frobenius error is
$$
\mathbb{E}\left(\lVert \mathsf{CUR}-M\rVert_F^2\right)
= \sum_{I\in \mathcal{S}_r(m)} \sum_{J\in \mathcal{S}_k(n)} p(I,J)\,\lVert M_{:,J} M_{I,J}^+ M_{I,:} - M\rVert_F^2.  
$$

The error analysis is simplified by introducing the submatrices $A,B,C,D$. 
We also use the following.
\begin{definition}
Let $I \in \mathcal{S}_r(m)$, $J \in \mathcal{S}_k(n)$ be selected using a sampling procedure (e.g., volume sampling).

Once $I$ and $J$ are chosen, we partition $M$ into blocks:
$$
M \sim \begin{bmatrix}
  A & B \\ C & D
\end{bmatrix}, \quad
A=M_{I,J}, B=M_{I,J^c}, C = M_{I^c,J}, D = M_{I^c,J^c}.$$
\end{definition}

The matrices $A,B,C,D$ correspond to the deterministic blocks introduced in Section~2.1 but now depend on the random index sets $I,J$ drawn from the sampling distribution. Each block is therefore a random matrix, and expectations below are taken with respect to this induced randomness.
Here the equivalence $\sim$ is understood as in the following definition
\begin{definition}[Equivalence under row and column permutations]
    Two matrices $M,\widetilde{M} \in \mathbb{R}^{m \times n}$ are said to be \emph{equivalent under permutations}, $$\widetilde{M}\sim M,$$ if there exist permutation matrices $P \in \mathbb{R}^{m \times m}$ and $Q \in \mathbb{R}^{n \times n}$ such that
$$
\widetilde{M}= P M Q.
$$
\end{definition}

\emph{Remarks.} 

The Frobenius norm is invariant under equivalence. For example, 
$$
\|\widetilde{M}\|_F = \|P M Q\|_F = \|M\|_F.
$$

The permutations $P, Q$ may be deterministic or random, depending on the context (e.g., when $\widetilde{M}$ is a random matrix depending on index sets $I$, $J$).

In particular, one has
$$
M \sim \begin{bmatrix} A & B \\ C & D \end{bmatrix}, \quad
\mathsf{R} \sim \begin{bmatrix} A & B \end{bmatrix}, \quad
\mathsf{C} \sim \begin{bmatrix} A \\ C \end{bmatrix}.
$$

\emph{Remark (on terminology).}
In this exposition we describe the matrices $A, B, C, D$ and the CUR approximation $\mathsf{CUR}$ as depending on sampled index sets $I, J$. We avoid calling them random matrices to maintain alignment with applied practice.

However, mathematically, since $I$ and $J$ are random (drawn from a sampling distribution such as volume sampling), the matrices constructed from them — including $A = M_{I,J}$, $B = M_{I,J^c}$, etc., and  $\mathsf{CUR}$ — are indeed random matrices. Our probabilistic statements (e.g., expectations over the error) are taken with respect to this induced randomness.
 
The following proposition expresses the deterministic CUR error in the block form of Section~2.3.
\begin{proposition}[Error formula for the CUR approximation]\label{prop:cur-error}
   Let $M \in \mathbb{R}^{m \times n}$, and let $\mathsf{CUR}$ be the CUR approximation as defined above. Then the squared Frobenius norm of the error is
$$
\|\mathsf{CUR}  - M\|_F^2 = \|(I - A A^+) B\|_F^2 + \|C A^+ B - D\|_F^2,
$$
valid wherever $A$ has full column rank.
\end{proposition}
\begin{proof}
    
Using the equivalences

$$
M \sim \begin{bmatrix} A & B \\ C & D \end{bmatrix}, \quad
\mathsf{C} \sim \begin{bmatrix} A \\ C \end{bmatrix}, \quad
\mathsf{R} \sim \begin{bmatrix} A & B \end{bmatrix},
$$

we compute

$$
\mathsf{C} \, \mathsf{U} \, \mathsf{R} \sim
\begin{bmatrix} A \\ C \end{bmatrix} A^+ \begin{bmatrix} A & B \end{bmatrix}
= \begin{bmatrix} A A^+ A & A A^+ B \\ C A^+ A & C A^+ B \end{bmatrix}.
$$

Subtracting the block form of $M$ gives

$$
\mathsf{CUR} - M \sim
\begin{bmatrix}
A A^+ A - A & A A^+ B - B \\
C A^+ A - C & C A^+ B - D
\end{bmatrix}.
$$

Since $A A^+ A = A$ and $C A^+ A = C$ (i.e.\ $\det(A^\top A) \neq 0$), we have

$$
\mathsf{C} \, \mathsf{U} \, \mathsf{R} - M \sim
\begin{bmatrix}
0 & (I - A A^+) B \\
0 & D - C A^+ B
\end{bmatrix}\quad \text{a.e},
$$

so the squared Frobenius norm is

$$
\|\mathsf{C} \, \mathsf{U} \, \mathsf{R} - M\|_F^2 = \|(I - A A^+) B\|_F^2 + \|D - C A^+ B \|_F^2 \quad \text{a.e.}
$$
\end{proof}

Hence the CUR error decomposes into orthogonal contributions from the column and row extensions, mirroring the local formula.

\subsection{Error of $B$ component using CUR approximation of $M$}

The following theorem extends the local column-addition identity to the expected global error of the $B$ block.
\begin{theorem}[Error of B component]\label{thm:error-B}
Let $M \in \mathbb{R}^{m \times n}$, $I \in \mathcal{S}_r(m)$ and $J \in \mathcal{S}_k(n)$, where $r \geq k$. Then
$$
\mathbb{E}\left( \| AA^+ B - B \|_F^2 \right) = \frac{(k+1)(r-k)}{m-k} \frac{\| C_{k+1}(M) \|_F^2}{\| C_k(M) \|_F^2}.
$$
\end{theorem}

\begin{proof}
From the definition of the expectation,

\begin{equation}
\mathbb{E}\left( \| AA^+ B - B \|_F^2 \right) 
 =\zeta^{-1}\sum_{I\in \mathcal{S}_r(m)} \sum_{J\in \mathcal{S}_k(n)} \vol(M_{I,J})^2\Vert(I-M_{I,J}M_{I,J}^+)M_{I,J^c}\Vert_F^2
\end{equation}

where $\zeta$ is the volume sampling normalisation factor and, on applying Proposition~\ref{prop:add-col} for each $j \in J^c$, one obtains

$$
\vol(M_{I,J})^2 \| (I-M_{I,J}M_{I,J}^+)M_{I,j} \|_2^2 = \vol(M_{I,K})^2,
$$

where $K = J \cup \{j\} \in \mathcal{S}_{k+1}(n)$. Then, on summing over all $j \in J^c$, one obtains

$$
\vol(M_{I,J})^2 \| (I-M_{I,J}M_{I,J}^+)M_{I,J^c} \|_2^2 
= \sum_{K \in \mathcal{S}_{k+1}(n)} \mathbf{1}_{J \subset K} \vol(M_{I,K})^2.
$$

Thus,
$$
\mathbb{E}\left( \| AA^+ B - B \|_F^2 \right) 
= \zeta^{-1} \sum_{J \in \mathcal{S}_k(n)} \sum_{K \in \mathcal{S}_{k+1}(n)} \sum_{I \in \mathcal{S}_r(m)} \mathbf{1}_{J \subset K} \vol(M_{I,K})^2.
$$

Now, applying the Cauchy--Binet formula,
$$
\vol(M_{I,K})^2 = \sum_{H \in \mathcal{S}_{k+1}(m)} \mathbf{1}_{H \subseteq I} \vol(M_{H,K})^2,
$$

and hence
\begin{equation*}
\begin{split}
\mathbb{E}\left( \| AA^+ B - B \|_F^2 \right)
&= \zeta^{-1} \sum_{J \in \mathcal{S}_k(n)} \sum_{K \in \mathcal{S}_{k+1}(n)} \\
&\quad\times \sum_{H \in \mathcal{S}_{k+1}(m)} \sum_{I \in \mathcal{S}_r(m)} \mathbf{1}_{J \subset K} \mathbf{1}_{H \subseteq I} \vol(M_{H,K})^2.
\end{split}
\end{equation*}

Reordering the summations (with $J$ outermost), one obtains
\begin{equation*}
\begin{split}
\mathbb{E}\left( \| AA^+ B - B \|_F^2 \right)
&= \zeta^{-1} \sum_{J \in \mathcal{S}_k(n)} \sum_{K \in \mathcal{S}_{k+1}(n)} \mathbf{1}_{J \subset K} \\
&\quad\times \sum_{H \in \mathcal{S}_{k+1}(m)} \left( \sum_{I \in \mathcal{S}_r(m)} \mathbf{1}_{H \subseteq I} \right) \vol(M_{H,K})^2.
\end{split}
\end{equation*}

The inner sum counts the number of sets $I$ of size $r$ containing $H$, and equals $\binom{m-k-1}{m-r}$. Thus,

\begin{equation*}
\begin{split}
\mathbb{E}\left( \| AA^+ B - B \|_F^2 \right)
&= \zeta^{-1} \binom{m-k-1}{m-r} \sum_{J \in \mathcal{S}_k(n)} \sum_{K \in \mathcal{S}_{k+1}(n)} \\
&\quad\times \mathbf{1}_{J \subset K} \sum_{H \in \mathcal{S}_{k+1}(m)} \vol(M_{H,K})^2.
\end{split}
\end{equation*}

Now, for each $K \in \mathcal{S}_{k+1}(n)$, there are exactly $k+1$ subsets $J \subset K$ of size $k$. Therefore,
$$
\sum_{J \in \mathcal{S}_k(n)} \mathbf{1}_{J \subset K} = k+1,
$$

and hence,
$$
\mathbb{E}\left( \| AA^+ B - B \|_F^2 \right) 
 = \zeta^{-1}(k+1) \binom{m-k-1}{m-r} \sum_{K \in \mathcal{S}_{k+1}(n)} \sum_{H \in \mathcal{S}_{k+1}(m)} \vol(M_{H,K})^2.
$$

Finally, it now follows from the Frobenius norm of the $(k+1)$-st compound matrix of $M$, that
$$
\| C_{k+1}(M) \|_F^2 = \sum_{H \in \mathcal{S}_{k+1}(m)} \sum_{K \in \mathcal{S}_{k+1}(n)} \vol(M_{H,K})^2,
$$
and, using Theorem~\ref{thm:volume-sampling-normalisation} for $\zeta$, we conclude that
$$
\mathbb{E}\left( \| AA^+ B - B \|_F^2 \right) = \frac{(k+1)(r-k)}{m-k} \frac{\| C_{k+1}(M) \|_F^2}{\| C_k(M) \|_F^2}.
$$
\end{proof}

The expectation of the $B$-component error thus inherits the same geometric structure as the local residual in Section~2.1.

The next theorem provides the companion estimate for the $D$ block.
\begin{theorem}[Error of approximating $D$]\label{thm:error-D}
Let $M \in \mathbb{R}^{m \times n}$, and let $I \in \mathcal{S}_r(m)$, $J \in \mathcal{S}_k(n)$ with $r \geq k$. Then
\begin{equation} \label{eq:errD}
  \mathbb{E}\left( \| CA^+B - D \|_F^2 \right) \leq \frac{(k+1)^2(m-r)}{m-k} \frac{\| C_{k+1}(M) \|_F^2}{\| C_k(M) \|_F^2}.
\end{equation}
\end{theorem}

\begin{proof}
Consider elements outside $I$ and $J$:

\begin{itemize}
\item Pick $i \in I^c$ (a row not in $I$),
\item Pick $j \in J^c$ (a column not in $J$).
\end{itemize}

Define

\begin{itemize}
\item $K = I \cup \{i\} \in \mathcal{S}_{r+1}(m)$,
\item $L = J \cup \{j\} \in \mathcal{S}_{k+1}(n)$.
\end{itemize}

Applying Proposition 3 yields

$$
\mathbb{E}\left( \| CA^+B - D \|_F^2 \right) 
\leq \zeta^{-1} \sum_{I \in \mathcal{S}_r(m)} \sum_{J\in \mathcal{S}_k(n)} \sum_{i \in I^c} \sum_{j \in J^c} \left( \vol(M_{K,L})^2 - \vol(M_{I,L})^2 \right).
$$

Proceed in two steps:

First fix $L \in \mathcal{S}_{k+1}(n)$ and expand the difference using the Cauchy--Binet formula:
$$
\det(M_{K,L}^\top M_{K,L}) = \sum_{H \in \mathcal{S}_{k+1}(m)} \mathbf{1}_{H \subseteq K} \vol(M_{H,L})^2,
$$
and
$$
\det(M_{I,L}^\top M_{I,L}) = \sum_{H \in \mathcal{S}_{k+1}(m)} \mathbf{1}_{H \subseteq I} \vol(M_{H,L})^2.
$$

Thus,
$$
\det(M_{K,L}^\top M_{K,L}) - \det(M_{I,L}^\top M_{I,L}) 
= \sum_{H \in \mathcal{S}_{k+1}(m)} \mathbf{1}_{H \subseteq K} \mathbf{1}_{H \not\subseteq I} \vol(M_{H,L})^2.
$$

Fix $H \in \mathcal{S}_{k+1}(m)$ and evaluate

$$
\sum_{K \in \mathcal{S}_{r+1}(m)} \mathbf{1}_{H \subseteq K} \sum_{\substack{I \subset K \\ |I|=r}}  \mathbf{1}_{H \not\subseteq I}.
$$

First sum over $I$ and then over $K$.
\begin{itemize}
\item For each fixed $K$, sum over $I \subset K$ with $|I|=r$.
Each $I$ corresponds to omitting exactly one element from $K$.

\item The subsets $I$ are in bijection with the elements $i \in K$ (the element omitted to form $I$).

\item We ask: when does $H \not\subseteq I$?
\begin{itemize}
\item If $i \in H$, then $H$ is missing an element and thus $H \not\subseteq I$.
\item If $i \notin H$, then $H \subseteq I$.
\end{itemize}

\item Thus, $i \in H$ $\quad \Longleftrightarrow \quad$ $H \not\subseteq I$.

\item Since $|H| = k+1$, there are exactly $k+1$ elements $i \in H$.

\item Thus, for each $K$, the number of subsets $I$ such that $H \not\subseteq I$ is $k+1$.
\end{itemize}

Thus,
$$
\sum_{\substack{I \subset K \\ |I|=r}} \mathbf{1}_{H \not\subseteq I} = k+1
$$
for each fixed $K$ with $H \subseteq K$.
Now sum over $K$ with $H\subseteq K$.
The number of such $K$ is $\binom{m-k-1}{m-r-1}$.

Hence, the total contribution over $K$ and $I$ is $(k+1) \binom{m-k-1}{m-r-1}$.

Now sum over $J$:

\begin{itemize}
\item For each fixed $L \in \mathcal{T}_{k+1}(n)$,
\item The number of subsets $J \subset L$ of size $k$ is $\binom{k+1}{k} = k+1$.
\end{itemize}

Thus, the contribution from summing over $J$ is another factor of $k+1$.

---

Putting everything together, we obtain
$$
\mathbb{E}\left( \| CA^+B - D \|_F^2 \right) \leq \frac{(k+1)^2}{\zeta} \binom{m-k-1}{m-r-1} 
\sum_{H \in \mathcal{S}_{k+1}(m)} \sum_{L \in \mathcal{S}_{k+1}(n)} \vol(M_{H,L})^2.
$$

Recognizing the Frobenius norm,
$$
\| C_{k+1}(M) \|_F^2 = \sum_{H \in \mathcal{S}_{k+1}(m)} \sum_{L \in \mathcal{S}_{k+1}(n)} \vol(M_{H,L})^2,
$$
and using Theorem 1 for $\zeta$,
$$
\zeta = \binom{m-k}{r-k} \| C_k(M) \|_F^2,
$$
we conclude
$$
\mathbb{E}\left( \| CA^+B - D \|_F^2 \right) \leq \frac{(k+1)^2(m-r)}{m-k} \frac{\| C_{k+1}(M) \|_F^2}{\| C_k(M) \|_F^2}.
$$
\end{proof}

This bound is the stochastic analogue of the symmetric case treated in Proposition~6 of Section~2.3.

\subsection{Error of the CUR approximation of $M$}

Combining Theorems~\ref{thm:error-B} and~\ref{thm:error-D} yields a complete interpolation bound for the expected CUR error.
\begin{theorem}[Interpolation error bound for oversampling]\label{thm:interpolation-bound}
Let $M \in \mathbb{R}^{m \times n}$ and $r \geq k$. Then
$$
\mathbb{E}\left( \lVert \mathsf{C}\mathsf{U}\mathsf{R} - M \rVert_F^2 \right)
\leq \left(\frac{m-r}{m-k}(k+1)^2 + \frac{r-k}{m-k}(k+1) \right) \frac{ \| C_{k+1}(M) \|_F^2 }{ \| C_k(M) \|_F^2 }.
$$
\end{theorem}

\begin{proof}
By Theorem ~\ref{thm:error-B}, we have
$$
\mathbb{E}\left( \| AA^+B - B \|_F^2 \right) = \frac{(k+1)(r-k)}{m-k} \frac{ \| C_{k+1}(M) \|_F^2 }{ \| C_k(M) \|_F^2 }.
$$
By Theorem~\ref{thm:error-D}, we have
$$
\mathbb{E}\left( \| CA^+B - D \|_F^2 \right) \leq \frac{(k+1)^2(m-r)}{m-k} \frac{ \| C_{k+1}(M) \|_F^2 }{ \| C_k(M) \|_F^2 }.
$$
Adding the two bounds gives
\begin{equation*}
\begin{split}
\mathbb{E}\left( \| CA^+B - D \|_F^2 + \| AA^+B - B \|_F^2 \right)
&\leq \frac{(k+1)}{m-k} \big( (r-k) + (k+1)(m-r) \big) \\
&\quad\times \frac{ \| C_{k+1}(M) \|_F^2 }{ \| C_k(M) \|_F^2 }.
\end{split}
\end{equation*}
The result then follows from Proposition~\ref{prop:cur-error}.
\end{proof}

As a function of $r$, this bound decreases linearly from $(k+1)^2$ at $r=k$ to $(k+1)$ at $r=m$.

\medskip
\noindent\textbf{Remark (Positive definite case).}
When $M$ is symmetric positive definite, Proposition~\ref{prop:rect-bound} in Section~2 shows that the two
error components coincide. In this case, the bound in Theorem \ref{thm:interpolation-bound}
becomes an equality: both the randomised CUR method (Nyström) and the general CUR approximation
discussed here have the same expected Frobenius error bound.
This unifies the symmetric and general cases under the same determinant framework.

\begin{proposition}[Error bounds in terms of singular values]\label{prop:singular-value-bounds}
Let $M\in \R^{m\times n}$ be a given matrix and $k\leq r < min(m,n)$. Then

\[\mathbb{E}\left(\|M-M_k\|_F^2\right)\leq
\left(\frac{m-r}{m-k}(k+1)^2 + \frac{r-k}{m-k}(k+1) \right)
\frac{e_{k+1}(\sigma_1^2,\cdots,\sigma_m^2)}{e_k(\sigma_1^2,\cdots,\sigma_m^2)}\]

where $\sigma_i$ is the $i$-th singular value of $M$. 
\end{proposition}

\begin{proof}
The proof is a direct application of Theorem~\ref{thm:interpolation-bound} (cf.\ Lemma~\ref{lem:gram-increment}) and  the fact that
$$\| C_k(M) \|_F^2 = e_k(\sigma_1^2,\cdots,\sigma_m^2).$$
\end{proof}

It then follows from Lemma 1 that:
\begin{corollary}
Let $M\in \R^{m,n}$ be a given matrix. Then

\[E\left(\|M-M_k\|_F^2\right)\leq
 \left(\frac{m-r}{m-k}(k+1)^2 + \frac{r-k}{m-k}(k+1) \right)
\sum_{i=k+1}^n \sigma_i^2\]
\end{corollary}

The corollary provides some idea of how increasing the $r$ may reduce the error. But there is
still this factor $k+1$ in the error bound for the Frobenius norm. 
This factor is known to be tight; see~\cite{deshpande2006matrix}.
\fussy

\printbibliography

\end{document}